\documentclass[11pt]{article}
\usepackage{graphicx, amsmath}
\usepackage[english]{babel}
\usepackage{amsmath}
\usepackage{graphicx}
\usepackage{amsfonts}
\usepackage{amssymb}

\newtheorem{theorem}{Theorem}
\newtheorem{corollary}[theorem]{Corollary}
\newtheorem{lemma}[theorem]{Lemma}
\newtheorem{proposition}[theorem]{Proposition}
\newtheorem{definition}[theorem]{Definition}
\newtheorem{remark}{Remark}
\newtheorem{example}{Example}

\newcommand{\R}{\mathbb{R}}
\begin{document}

%%%%%%%%%%%%%%%%%%%%%%%%%%%%%%%%%%%%%%%%%%%%%%%%%%%%%%%%%%%%
%%%%                                                    %%%%
%%%%     Clarke subgradients of stratifiable functions  %%%%
%%%%   J. Bolte, A. Daniilidis, A. Lewis, M. Shiota     %%%%
%%%%                                                    %%%%
%%%%               January 18, 2006                     %%%%
%%%%            [Last version arisd]                    %%%%
%%%%                                                    %%%%
%%%%%%%%%%%%%%%%%%%%%%%%%%%%%%%%%%%%%%%%%%%%%%%%%%%%%%%%%%%%

\begin{center}
{\LARGE  Clarke subgradients of stratifiable functions}

\vspace{0.2cm}

\vspace{0.8cm}\textbf{J\'{e}r\^{o}me BOLTE, Aris DANIILIDIS, Adrian LEWIS \&
Masahiro SHIOTA}

\bigskip
\end{center}

\noindent\textbf{Abstract }We establish the following result: if
the graph of a (nonsmooth) real-extended-valued function
$f:\mathbb{R}^{n}\rightarrow \mathbb{R}\cup\{+\infty\}$ is closed
and admits a Whitney stratification, then the norm of the
 gradient of $f$ at $x\in\mbox{dom\,}f$ relative to the
stratum containing $x$ bounds from below all norms of Clarke
subgradients of $f$ at $x$. As a consequence, we obtain some
Morse-Sard type theorems as well as a nonsmooth Kurdyka-\L
ojasiewicz inequality for functions definable in an arbitrary
o-minimal structure.

\bigskip


\bigskip

\textit{Secondary} 26D10, 32B20

\section{Introduction}

\hspace{0.4cm}Nonsmoothness in optimization seldom occurs in an arbitrary
manner, but instead it is often well-structured in the sense that a naturally
arising manifold $\mathcal{M}$ contains the minimizer, and the function is
smooth along this manifold. We quote \cite{Lewis2002} for formal definitions,
examples and more details. In the last two decades, several researchers have
tried to capture this intuitive idea in order to develop algorithms ensuring
better convergence results: see for instance the pioneer work \cite{LOS2000},
and also \cite{MM2004} for a recent survey.

\smallskip

In this work we shall be interested in a particular class of
well-structured (nonsmooth) functions, namely functions admitting
a Whitney stratification (see Section 2 for definitions). Since
this class contains in particular the semialgebraic and the
subanalytic functions (more generally, functions that are
definable in some o-minimal structure over $\R$), the derived results can
directly be applied in several concrete optimization problems
involving such structures. Our central idea is to relate
derivative ideas from two distinct mathematical sources:
variational analysis and differential geometry. Specifically, we
derive a lower bound on the norms of Clarke subgradients at a
given point in terms of the ``Riemannian'' gradient with respect to
the stratum containing that point. This is a direct consequence of
the ``projection formula'' given in Proposition \ref{projec} and
has as corollaries a Morse-Sard type theorem for Clarke critical
points of lower semicontinuous Whitney stratifiable functions
(Corollary \ref{cor1}(ii)) as well as a nonsmooth version of the
Kurdyka-\L ojasiewicz inequality for lower semicontinuous
definable functions (Theorem \ref{Theorem_KLineq}). Although the
proofs are reasonably routine, analogous results fail for the
(broader) convex-stable subdifferential (introduced and studied in
\cite{BLO2002}), unless $f$ is assumed to be locally Lipschitz
continuous, see Remark \ref{Remark_convex-stable}.

\smallskip

As mentioned above, convergence theory for algorithms is one of the main
motivations for this work. In order to treat nonconvex (and nonsmooth)
minimization problems, the authors of \cite{BLO2002} introduced an algorithm
called the ``gradient sampling algorithm''. The idea behind this algorithm was
to sample gradients of nearby points of the current iterate and to produce the
next iterate by following the vector of minimum norm in the convex hull
generated by the sampled negative gradients. In the case that the function is
locally Lipschitz, the above method can be viewed as a kind of $\varepsilon$-Clarke
subgradient algorithm for which both theoretical and numerical results are
quite satisfactory, see \cite{BLO2002}. The convergence of the whole sequence
of iterates remains however an open question and this is also
 the case for many classical subgradient methods for nonconvex minimization,
 see \cite{Kiwiel}. Following the ideas of
\cite{Loja63} and \cite{Kurdyka98}, we would hope that \L ojasiewicz
inequality that we develop ((\ref{KL}) in Section \ref{Kurdik}) could play a
prominent role in the global convergence of subgradient algorithms.

\section{Preliminaries}

\hspace{0.40cm} In this section we recall several definitions and results
concerning nonsmooth analysis (subgradients, generalized critical points) and
stratification theory. In what follows the vector space $\mathbb{R}^{n}$ is
endowed with its canonical scalar product $\langle\cdot,\cdot\rangle$.

\bigskip

\textbf{Nonsmooth analysis.} Given an extended-real-valued function
$f:\mathbb{R}^{n}\rightarrow\mathbb{R}\cup\{+\infty\}$ we denote its domain by
$\mbox{dom\,}f:=\left\{  x\in\mathbb{R}^{n}:f(x)<+\infty\right\}  $, its graph
by
\[
\mbox{Graph\,}f:=\left\{  (x,f(x))\in\mathbb{R}^{n}\mathbb{\times R}:x\in\mbox
{dom\,}f\right\}
\]
and its epigraph by
\[
\mbox{epi\,}f:=\{(x,\beta)\in\mathbb{R}^{n}\mathbb{\times R}: f(x)\leq\beta\}.
\]
In this work we shall deal with \emph{lower semicontinuous}
functions, that is, functions for which $\mbox{epi\,}f$ is a
closed subset of $\mathbb{R} ^{n}\mathbb{\times R}$. In this
setting, we say that $x^{\ast}\in \mathbb{R}^{n}$ is a
\emph{Fr\'{e}chet subgradient} of $f$ at $x\in \mbox{dom\,}f$ provided
that
\begin{equation}
\underset{y\rightarrow x,y\neq
x}{\lim\inf}\;\frac{f(y)-f(x)-\langle x^{\ast },y-x\rangle}{\Vert
y-x\Vert}\text{ }\geq\text{ }0\text{.}\label{Frechet}
\end{equation}
The set of all Fr\'{e}chet subgradients of $f$ at $x$ is called the
\emph{Fr\'{e}chet subdifferential} of $f$ at $x$ and is denoted by
$\hat{\partial}f(x)$. If $x\notin\mathrm{dom\,}f$ then we set $\hat{\partial
}f(x)=\emptyset$.

\smallskip

Let us give a geometrical interpretation of the above definition:
it is well known that the gradient of a $C^{1}$ function
$f:\mathbb{R}^{n}\rightarrow \mathbb{R}$ at $x\in\mathbb{R}^{n}$
can be defined geometrically as the vector $\nabla
f(x)\in\mathbb{R}^{n}$ such that $(\nabla f(x),-1)$ is normal to
the tangent space $\mathrm{T}_{(x,f(x))}\mathrm{Graph\,}f$ of (the
$C^{1}$ manifold) $\mbox{Graph\,}f$ at $(x,f(x))$, that is,
\[
(\nabla f(x),-1)\;{\LARGE \bot}\;\mathrm{T}_{(x,f(x))}\mathrm{Graph\,}f.
\]
A similar interpretation can be stated for Fr\'{e}chet subgradients. Let us
first define the \emph{(Fr\'{e}chet) normal cone }of a subset $C$ of
$\mathbb{R}^{n}$ at $x\in C$ by
\begin{equation}
\hat{N}_{C}(x)=\left\{  v\in\mathbb{R}^{n}:\;\underset{y\in
C\setminus \{x\}}{\underset{y\rightarrow
x}{\limsup}}\,\langle\,v,\frac{y-x}
{||x-y||}\,\rangle\,\leq\,0\right\}  . \label{normalcone}
\end{equation}
Then it can be proved (see \cite[Theorem 8.9]{Rock98}, for example) that for a
nonsmooth function $f$ we have:
\begin{equation}
x^{\ast}\in\hat{\partial}f(x)\quad\mbox{if and only
if}\quad(x^{\ast},-1)\in\hat{N}_{\mathrm{epi\,}f}(x,f(x)).
\label{F}
\end{equation}
The Fr\'{e}chet subdifferential extends the notion of a derivative in the
sense that if $f$ is differentiable at $x$ then $\hat{\partial}f(x)=\{\nabla
f(x)\}$. However, it is not completely satisfactory in optimization, since
$\hat{\partial}f(x)$ might be empty-valued at points of particular interest
(think of the example of the function $f(x)=-||x||,$ at $x=0$). Moreover, the
Fr\'{e}chet subdifferential is not a closed mapping, so it is unstable
computationally. For this reason we also consider (see \cite[Chapter
8]{Rock98}, for example):

\begin{itemize}
\item [(i)]the \emph{limiting} subdifferential $\partial f(x)$ of $f$ at
$x\in\mbox{dom\,}f$:
\begin{equation}
p\in\partial f(x)\;\Longleftrightarrow\;\exists(x_{n},x_{n}^{\ast}
)_{n\in\mathbb{N}}\subset\mbox{Graph\,}{\hat{\partial}}f\,:\,\,\left\{
\begin{array}
[c]{l}
\underset{n\rightarrow\infty}{\lim}x_{n}=x,\\
\\
\underset{n\rightarrow\infty}{\lim}f(x_{n})=f(x),\\
\\
\underset{n\rightarrow\infty}{\lim}x_{n}^{\ast}=p,
\end{array}
\right.  \,\label{limiting}
\end{equation}
where $\mbox{Graph\,}{\hat{\partial}}f:=\{(u,u^*):u\in \hat{\partial}f(u)\}.$

\item[(ii)] the \emph{asymptotic limiting }subdifferential $\partial^{\infty
}f(x)$ of $f$ at $x\in\mbox{dom\,}f$:
\begin{equation}
q\in\partial^{\infty}f(x)\;\Longleftrightarrow\;\exists(y_{n},y_{n}^{\ast
})_{n\in\mathbb{N}}\subset\mbox{Graph\,}{\hat{\partial}}f,\,\exists
t_{n}\searrow0^{+}:\,\left\{
\begin{array}
[c]{l}
\underset{n\rightarrow\infty}{\lim}y_{n}=x,\\
\\
\underset{n\rightarrow\infty}{\lim}f(y_{n})=f(x),\\
\\
\underset{n\rightarrow\infty}{\lim}t_{n}y_{n}^{\ast}=q.
\end{array}
\right.  \label{limiting-asympt}
\end{equation}
\end{itemize}

When $x\notin\mbox{dom\,}f$ we set $\partial f(x)=\partial^{\infty
}f(x)=\emptyset$.\smallskip

The \emph{Clarke subdifferential} $\partial^{\circ}f(x)$ of $f$ at $x\in
\mbox{dom\,}f$ is the central notion of this work. It can be defined in
several (equivalent) ways, see \cite{Clarke1983}. The definition below (see
\cite[Theorem 8.11]{MordukShao96}) is the most convenient for our purposes.
(For any subset $S$ of $\mathbb{R}^{n}$ we denote by $\overline{\mbox{co}}\,S$
the closed convex hull of $S.$)

\begin{definition}
[Clarke subdifferential]{\label{Definition_Clarke}} \rm{The Clarke
subdifferential $\partial^{\circ}f(x)$ of $f$ at $x$ is the set}
\begin{equation}
\partial^{\circ}f(x)=\left\{
\begin{array}
[c]{cc} \overline{\mbox{co}}\,\left\{  \partial
f(x)+\partial^{\infty}f(x)\right\}
, & \quad\text{if }\,x\,\in\,\mathrm{\mbox{dom\,}}f\\
& \\
\mathrm{{{\emptyset,}}} & \quad\text{if }\,x\,\notin\,\mathrm{\mbox{dom\,}}f
\end{array}
\right.  \label{clarke}
\end{equation}
\end{definition}

\begin{remark}
\label{Rem_adrian}\rm{It can be shown that an analogous to
(\ref{F}) formula holds also for the Clarke subdifferential, if
$\hat{N}_{\mathrm{epi\,}f}(x,f(x))$ is replaced by the Clarke
normal cone, which is the closed convex hull of the limiting
normal cone. The latter cone comes naturally from the Fr\'{e}chet
normal cone by closing its graph, see \cite[pp. 305, 336]{Rock98}
for details.}
\end{remark}

{F}rom the above definitions it follows directly that for all
$x\in \mathbb{R}^{n}$, one has
\begin{equation}
\hat{\partial}f(x)\subset\partial f(x)\subset\partial^{\circ}f(x).
\label{inclusions}
\end{equation}
The elements of the limiting (respectively, Clarke) subdifferential are called
limiting (respectively, Clarke) subgradients.

\smallskip

The notion of a Clarke critical point (respectively, critical value,
asymptotic critical value) is defined as follows.

\begin{definition}
[Clarke critical point]\label{Crit} \rm{We say that
$x\in\mathbb{R}^{n}$ is a \emph{Clarke critical point} of the
function $f$ if}
\[
\partial^{\circ}f(x)\ni0.
\]
\end{definition}

\begin{definition}
[(asymptotic) Clarke critical value]\label{def-crit-value} \rm{(i)
We say that $\alpha\in\mathbb{R}$ is a \emph{Clarke critical
value} of $f$ if the level set $f^{-1}(\{\alpha\})$ contains a
Clarke critical point.} \smallskip

\noindent\rm{(ii) We say that
$\lambda\in\mathbb{R}\cup\{\pm\infty\}$ is an \emph{asymptotic
Clarke critical value} of $f$, if there exists a sequence
$(x_{n},x_{n}^{\ast})_{n\geq1}\subset\mbox{Graph\,}{\partial^{\circ}}f$,
such that}
\[
\left\{
\begin{array}
[c]{c}
f(x_{n})\,\rightarrow\,\lambda\\
\\
(1\,+\,||x_{n}||)\,||x_{n}^{\ast}||\,\rightarrow\,0.
\end{array}
\right.
\]
\end{definition}

Let us make some observations concerning the above definitions:

\begin{remark}
\label{Remark1} \rm{(i) Both limiting and Clarke subgradients are
generalizations of the usual gradients: indeed, if $f$ is $C^{1}$
around $x$ then we have:}
\[
\partial f(x)=\partial^{\circ}f(x)=\{\nabla f(x)\}.
\]

\rm{(ii) The asymptotic limiting subdifferential should not be
thought as a set of subgradients. Roughly speaking it is designed
to detect ``vertical tangents'' to the graph of $f$. For instance,
for the (nonsmooth) function $f(x)=x^{\frac {1}{3}}$
($x\in\mathbb{R}$) we have $\partial^{\infty}f(0)=\mathbb{R}_{+}$.
Note that since the domain of the Fr\'{e}chet subdifferential is dense
in $\mathrm{dom\,}f$, we always have $\partial^{\infty}f(x)\ni0$,
for all $x\in\mathrm{dom\,}f$ (see also \cite[Corollary
8.10]{Rock98}); therefore, this latter relation cannot be regarded
as a meaningful definition of critical points.}

\smallskip

\rm{(iii) To illustrate the definition of the Clarke critical
point (Definition \ref{Definition_Clarke}) let us consider the
example of the function $f:\mathbb{R}\rightarrow\mathbb{R}$
defined by}
\[
f(x)=\left\{
\begin{array}
[c]{ll}
\quad x, & \text{if }x\leq0\\
-\sqrt{x} & \text{if }x>0.
\end{array}
\right.
\]
\rm{Then $\hat{\partial}f(0)=\emptyset$ and $\partial f(0)=\{1\}$.
However, since $\partial^{\infty}f(0)=\mathbb{R}_{-}$ it follows
from (\ref{clarke}) that $\partial^{\circ} f(0)=(-\infty,1]$, so
$x=0$ is a Clarke critical point.}

\smallskip

\rm{(iv) It follows from Definition \ref{def-crit-value} that
every Clarke critical value $\alpha\in\mathbb{R}$ is also an
asymptotic Clarke critical value (indeed, given $x_{0}\in
f^{-1}(\{\alpha\})$ with $0\in
\partial^{\circ}f(x_{0})$, it is sufficient to take $x_{n}:=x_{0}$ and
$x_{n}^{\ast}=0$). Note that in case that $f$ has a bounded domain
$\mbox{dom\,}f$, Definition \ref{def-crit-value}\thinspace(ii) can
be simplified in the following way: the value
$\lambda\in\mathbb{R}\cup \{\pm\infty\}$ is critical if, and only
if, there exists a sequence
$(x_{n},x_{n}^{\ast})_{n\geq1}\subset\mbox{Graph\,}{\partial^{\circ}}f$,
such that $f(x_{n})\,\rightarrow\,\lambda$ and
$x_{n}^{\ast}\,\rightarrow\,0$.}
\end{remark}

\bigskip

\textbf{Stratification results. }By the term \emph{stratification} we mean a
locally finite partition of a given set into differentiable manifolds, which,
roughly speaking, fit together in a regular manner. Let us give a formal
definition of a $C^{p}$-stratification of a set (for general facts about
stratifications we quote \cite{Mather1970} or \cite{Kaloshin2005} and
references therein).

\smallskip

Let $X$ be a subset of $\mathbb{R}^{n}$ and $p$ a positive integer. A $C^{p}$
\emph{stratification} $\mathcal{X}=(X_{i})_{i\in I}$ of $X$ is a locally
finite partition of $X$ into $C^{p}$ submanifolds $X_{i}$ of $\mathbb{R}^{n}$
such that for each $i\neq j$
\[
\overline{X_{i}}\cap X_{j}\neq\emptyset\,\Longrightarrow\,X_{j}\subset
\overline{X_{i}}\setminus X_{i}.
\]
The submanifolds $X_{i}$ are called \emph{strata} of $\mathcal{X}$.
Furthermore, given a finite collection $\{A_{1},\ldots,A_{q}\}$ of subsets of
$X$, a stratification $\mathcal{X=}(X_{i})_{i\in I}$ is said to be
\emph{compatible with the collection} $\{A_{1},\ldots,A_{q}\}$ if each $A_{i}$
is a locally finite union of strata $X_{j}$.

\smallskip

In this work we shall use a special type of stratifications
(called Whitney stratifications) for which the strata are such
that their tangent spaces also ``fit regularly''. To give a
precise meaning to this statement, let us first define the
\emph{distance} (or $\emph{gap}$) of two vector subspaces $V$ and
$W$ of $\mathbb{R}^{n}$ by the following standard formula
\[
D(V,W)=\mathrm{max\;}\left\{  \mathrm{\,}\sup_{v\in V,\,||v||=1}
d(v,W),\mathrm{\;}\sup_{w\in W,\,||w||=1}d(w,V)\right\}  .
\]
Note that
\[
\sup_{v\in V,\,||v||=1}d(v,W)=0\Longleftrightarrow V\subset W.
\]
Further we say that a sequence $\{V_{k}\}_{k\in\mathbb{N}}$ of subspaces of
$\mathbb{R}^{n}$ \emph{converges} to the subspace $V$ of $\mathbb{R}^{n}$ (in
short, $V=\underset{k\rightarrow+\infty}{\lim}V_{k}$)\ provided
\[
\underset{k\rightarrow+\infty}{\lim}D(V_{k},V)=0.
\]
Notice that in this case all the subspaces $V_{k}$ eventually have the same
dimension (say $d$), so that the above convergence is essentially equivalent
to the convergence in the grassmannian manifold $G_{d}^{n}$.

\smallskip

A $C^{p}$-stratification $\mathcal{X}=(X_{i})_{i\in I}$ of $X$ has the
\emph{Whitney-}($a$)\emph{ property,} if for each $x\in\overline{X_{i}}\cap
X_{j}$ (with $i\neq j$) and for each sequence $\{x_{k}\}_{k\geq1}\subset
X_{i}$ we have:
\[
\left.
\begin{array}
[c]{ll}
& \underset{k\rightarrow\infty}{\lim}\mathcal{\;}x_{k}\mathcal{\;}=x\\
\text{and} & \\
&
\underset{k\rightarrow\infty}{\lim}\mathcal{\;}T_{x_{k}}X_{i}\mathcal{\;}
=\mathcal{T}
\end{array}
\right\}  \mathcal{\;}\Longrightarrow\mathcal{\;}T_{x}X_{j}\mathcal{\;}
\subset\mathcal{\;T}
\]
where $T_{x}X_{j}$ (respectively, $T_{x_{k}}X_{i}$) denotes the tangent space
of the manifold $X_{j}$ at $x$ (respectively, of $X_{i}$ at $x_{k}$). In the
sequel we shall use the term \emph{Whitney stratification} to refer to a
$C^{1}$-stratification with the Whitney-($a$) property.

\section{Projections formulae for subgradients}

\label{section3}

\hspace{0.40cm} In this section we make precise the links between
the Clarke subgradients of a function whose graph (is closed and)
admits a Whitney stratification and the  gradients of
$f$ (with respect to the strata). As a corollary we obtain a
nonsmooth extension of the Morse-Sard theorem for such functions
(see Corollary \ref{cor1}).

\smallskip

Let $f:\mathbb{R}^{n}\rightarrow\mathbb{R} \cup\{+\infty\}$ be an
extended-real-valued function with a nonempty closed domain
$\mathrm{dom\,}f$ (that is, $f$ is lower semicontinuous). We shall
deal with Whitney stratifications $\mathcal{S}=(S_{i})_{i\in I}$
of the graph $\mbox{Graph\,}f$ of $f$ satisfying for all $i\in I$
and $u\in S_{i}$ the transversality condition:
\[
e_{n+1}\,\notin\,T_{u}S_{i}\hspace{2.4cm}\text{(}\mathcal{H}\text{)}
\]
where $$e_{n+1}=(0,\ldots,0,1)\in\mathbb{R}^{n+1}.$$

\begin{remark}
\label{Remark-adrian} \rm{If $f$ is locally Lipschitz continuous,
then it is easy to check that any stratification of  $\mbox{Graph\,}f$ 
must  automatically
satisfy ($\mathcal{H}$). This might also happen for other
functions (think of the non-locally Lipschitz function
$f(x)=\sqrt{|x|}$: every stratification of $\mbox{Graph\,}f$
should contain the stratum $S_{i}=\{(0,0)\}$), however the example
of the function $f(x)=x^{3}$ shows that this is not the case for
any (continuous stratifiable) function $f$ and any stratification
of its graph (consider the trivial stratification consisting of
the single stratum $\mbox{Graph\,}f$ and take $u=(0,0)$).}
\end{remark}

Let us denote by $\Pi:\mathbb{R}^{n+1}\rightarrow\mathbb{R}^{n}$
the canonical projection on $\mathbb{R}^{n}$, that is,
\[
\Pi(x_{1},\ldots,x_{n},t)=(x_{1},\ldots,x_{n}).
\]
For each $i\in I$ we set
\begin{equation}
X_{i}=\Pi(S_{i})\qquad\mbox{ and }\qquad
f_{i}=f|_{X_{i}}.\label{ww}
\end{equation}
Due to the above assumptions one has for all $i\in I$\smallskip

- $X_{i}$ is a $C^{1}$ submanifold of $\mathbb{R}^{n}$,\smallskip

- $f_{i}:X_{i}\rightarrow\mathbb{R}$ is a $C^{1}$ function.\smallskip

- $\mathcal{X}=(X_{i})_{i\in I}$ is a Whitney stratification of
$\mathrm{dom\,}f=\Pi(\mbox{Graph\,}f).$

\bigskip

\noindent\textbf{Notation.} In the sequel, for any $x\in\mbox{dom\,}f,$ we
shall denote by $X_{x}$ (respectively, $S_{x}$) the stratum of $\mathcal{X}$
(respectively of $\mathcal{S}$) containing $x$ (respectively $(x,f(x))$). The
manifolds $X_{i}$ are here endowed with the metric induced by the canonical
Euclidean scalar product of $\mathbb{R}^{n}$. Using the inherited Riemannian
structure of each stratum $X_{i}$ of $\mathcal{X}$, for any $x\in X_{i}$, we
denote by $\nabla_{\emph{R}}f(x)$ the  gradient of $f_{i}$ at
$x$ with respect to the stratum $X_{i},\langle\cdot,\cdot\rangle$.

\begin{proposition}
[Projection formula]\label{projec} Let
$f:\mathbb{R}^{n}\rightarrow \mathbb{R}\cup\{+\infty\}$ be a lower
semicontinuous function and assume that $\mbox{Graph\,}f$ admits a
Whitney stratification $\mathcal{S}=(S_{i})_{i\in I}$ satisfying
($\mathcal{H}$). Then for all $x\in\mathrm{dom\,}f$ we have
\begin{equation}
\mathrm{Proj\,}_{T_{x}X_{x}}\,\partial
f(x)\,\subset\,\{\nabla_{\emph{R}
}f(x)\}\quad;\quad\mathrm{Proj\,}_{T_{x}X_{x}}\,\partial^{\infty
}f(x)\,=\,\{0\}\label{proj-lim}
\end{equation}
and
\begin{equation}
\mathrm{Proj\,}_{T_{x}X_{x}}\mathrm{\,}\partial^{\circ}f(x)\mathrm{\,}
\subset\mathrm{\,}\{\nabla_{\emph{R}}f(x)\},\label{proj}
\end{equation}
where $\mathrm{Proj\,}_{\mathcal{V}}:\mathrm{\,}\mathbb{R}^{n}\rightarrow
\mathcal{V}$ denotes the orthogonal projection on the vector subspace
$\mathcal{V}$ of $\mathbb{R}^{n}$.
\end{proposition}

\noindent\textbf{Proof} We shall use the above notation (and in particular the
notation of (\ref{ww})).\newline Let us first describe the links between the
Fr\'{e}chet subdifferential $\hat{\partial}f(x)$ and the  gradient
of $f\,|_{X_{x}}$ at a point $x\in\mathrm{dom\,}f$. For any $v\in T_{x}X_{x}$
and any continuously differentiable curve $c:(-\varepsilon,\varepsilon
)\rightarrow X_{x}$ $(\varepsilon>0)$ with $c(0)=x$ and $\dot{c}(0)=v$, the
function
\[
f\circ c\,(:=f_{i}\circ
c)\,:(-\varepsilon,\varepsilon)\rightarrow\mathbb{R}
\]
is continuously differentiable. In view of \cite[Theorem 10.6, page
427]{Rock98}, we have
\[
\left\{  \,\langle x^{\ast},\,v\rangle:\,x^{\ast}\in\hat{\partial
}f(x)\,\right\}  \,\subset\,\left\{
\,\frac{d}{dt}\,f(c(t))\,|_{t=0} \,\right\}  .
\]
Since
$\frac{d}{dt}f(c(t))|_{t=0}=\langle\nabla_{\emph{R}}f(x),v\rangle$
it follows that
\begin{equation}
\mathrm{Proj\,}_{T_{x}X_{x}}\,\hat{\partial}f(x)\,\subset\,\{\nabla_{\emph{R}
}f(x)\}.\label{f1}
\end{equation}
In a second stage we prove successively that
\begin{equation}
\mathrm{Proj\,}_{T_{x}X_{x}}\,\partial
f(x)\,\subset\,\{\nabla_{\emph{R} }f(x)\}\text{\quad
and\quad}\mathrm{Proj\,}_{T_{x}X_{x}}\,\partial^{\infty
}f(x)\,\subset\,\{0\}.\label{l1}
\end{equation}
To this end, take $p\in\partial f(x)$, and let
$\{x_{k}\}_{k}\mathrm{\,} \subset\mathrm{\,dom\,}\hat{\partial}f$,
$x_{k}^{\ast}\in\hat{\partial} f(x_{k})$ be such that
$(x_{k},\,f(x_{k}))\rightarrow (x,\,f(x))$ and
$x_{k}^{\ast}\rightarrow p$. Due to the local finiteness property
of $\mathcal{S}$, we may suppose that the sequence
$\{u_{k}:=(x_{k} ,f(x_{k}))\}_{k}$ lies entirely in some stratum
$S_{i}$ of dimension $d$.\smallskip

If $S_{i}=S_{x}$ then by (\ref{f1}) we deduce that $\mathrm{Proj\,}
_{T_{x}X_{x}}\,(x_{k}^{\ast})\,=\nabla_{\emph{R}}f(x_{k})$, thus using the
continuity of the projection and the fact that $f\,|_{X_{x}}$ is $C^{1}$ (that
is, $\nabla_{\emph{R}}f(x_{k})\rightarrow\nabla_{\emph{R}}f(x)$) we obtain
$\mathrm{Proj\,}_{T_{x}X_{x}}\,(p)\,=\nabla_{\emph{R}}f(x)$.\smallskip

If $S_{i}\neq S_{x}$, then from the convergence
$(x_{k},\,f(x_{k}))\rightarrow (x,\,f(x))$ we deduce that
$\overline{S_{i}}\cap S_{x}\neq\emptyset$ (thus $d=\dim S_{i}>\dim
S_{x}$). Using the compactness of the grassmannian manifold
$G_{d}^{n}$, we may assume that the sequence
$\{T_{u_{k}}S_{i}\}_{k\geq1}$ converges to some vector space
$\mathcal{T}$ of dimension $d$. Then the Whitney-($a$) property
yields that $\mathcal{T}\supset T_{(x,f(x))}S_{x}$. Recalling
(\ref{F}), for each $k\geq1$ we have that the vector
$(x_{k}^{\ast},-1)$ is Fr\'{e}chet normal to the epigraph $\mbox
{epi\,}f$ of $f$ at $u_{k}$, hence it is also normal (in the
classical sense) to the tangent space $T_{u_{k}}S_{i}$. By a
standard continuity argument the vector
\[
(p,-1)=\,\underset{k\rightarrow\infty}{\lim}\,(x_{k}^{\ast},-1)
\]
must be normal to $\mathcal{T}$ and a fortiori to
$T_{(x,f(x))}S_{x}$. By projecting $(p,-1)$ orthogonally on
$T_{x}X_{x}\mathbb{\,}+\mathbb{\,R\,}
e_{n+1}\mathbb{\,}\supset\mathbb{\,}T_{(x,f(x))}S_{x}$, we notice
that $(\mathrm{Proj\,}_{T_{x}X_{x}}(p),-1)$ is still normal to
$T_{(x,f(x))}S_{x}$. By the definition of the 
subgradient we conclude that
\begin{equation}
\mathrm{Proj\,}_{T_{x}X_{x}}\,(p)=\nabla_{\emph{R}}f(x),\label{morduk}
\end{equation}
thus the first part of (\ref{l1}) follows.\smallskip

Let now any $q\in\partial^{\infty}f(x)$. By definition there exist
$\{y_{k}\}_{k}\mathrm{\,}\subset\mathrm{\,dom\,}\hat{\partial}f$,
$y_{k} ^{\ast}\in\hat{\partial}f(y_{k})$ and a positive sequence
$t_{k}\searrow0^{+}$ such that $(y_{k},\,f(y_{k}))\rightarrow
(y,\,f(y))$ and $t_{k} y_{k}^{\ast}\rightarrow q$. As above we may
assume that the sequence $\{y_{k}\}_{k}$ belongs to some stratum
$S_{i}$ and that the tangent spaces
$T_{u_{k}}S_{i}=T_{(x_{k},f(x_{k}))}S_{i}$ converge to some
$\mathcal{T}$. Since $t_{k}(y_{k}^{\ast},-1)$ is normal to
$T_{u_{k}}S_{i}$ we can similarly deduce that
$(\mathrm{Proj\,}_{T_{x}X_{x}}\,(q),\,0)$ is normal to
$T_{(x,f(x))}S_{x}$. Since $\mathrm{Proj\,}_{\mathbb{R}^{n}\times
\{0\}}T_{(x,f(x))}S_{x}=T_{x}X_{x}$ this implies that
$\partial^{\infty }f(x)\subset\left(  T_{x}X_{x}\right)  ^{\perp}$
and the second part of (\ref{l1}) is proved. It now follows from
(\ref{l1}) and Remark\thinspace \ref{Remark1}\thinspace(ii) that
(\ref{proj-lim}) holds.\smallskip

In order to conclude let us recall (Definition \ref{Definition_Clarke}) that
$\partial^{\circ}f(x)=\overline{\mbox{co}}\,(\partial f(x)+\partial^{\infty
}f(x))$. In view of (\ref{l1}) any element of $\mbox{co
}(\partial f(x)+\partial^{\infty}f(x))$ admits $\nabla_{\emph{R}}f(x)$ as a
projection onto $T_{x}X_{x}$. By taking the closure of the previous set we
obtain (\ref{proj}).$\hfill\Box$

\bigskip

\begin{remark}
\label{Remark_clarke_empty} \rm{The inclusion in (\ref{proj}) may
be strict (think of the function $f(x)=-||x||^{1/2}$ at $x=0$
where $\partial^{\circ}f(0)=\emptyset$). Of course, whenever
$\partial^{\circ}f(x)$ is nonempty (for example, if $f$ is locally
Lipschitz), under the assumptions of Proposition \ref{projec} we
have}
\[
\mathrm{Proj\,}_{T_{x}X_{x}}\partial^{\circ}f(x)=\{\nabla_{\emph{R}}f(x)\}.
\]
\end{remark}

\begin{corollary}
\label{cor1} Assume that the graph of $f$ is closed and admits a
$C^{p} $-Whitney stratification satisfying ($\mathcal{H}$). Then:
\smallskip

\noindent(i) for all $x\in\mbox{dom\,}\partial^{\circ}f$ we have
\begin{equation}
||\nabla_{\emph{R}}f(x)||\;\leq\;||x^{\ast}||,\qquad
\mathrm{\mathit{for\ all\
}}x^{\ast}\in\partial^{\circ}f(x).\label{proj-1}
\end{equation}

\noindent(ii) {\rm(\textbf{Morse-Sard theorem})} If $p\,\geq\,n$,
then the set of Clarke critical values of $f$ has Lebesgue measure
0.
\end{corollary}

\noindent\textbf{Proof} Assertion (i) is a direct consequence of
(\ref{proj}) of Proposition \ref{projec}. To prove (ii), set
$C:=[\partial^{\circ} f]^{-1}(\{0\} )=\{x\in\mathbb{R}^{n}:
\partial^{\circ} f(x)\ni0\}$. Since the set of strata is at most
countable, the restrictions of $f$ to each of those yield a
countable family $\{f_{n}\}_{n\in\mathbb{N}}$ of $C^{p}$
functions. In view of (\ref{proj-1}), we have that $C
\subset\cup_{n\in\mathbb{N}}\nabla f_{n}^{-1}(0)$. The result
follows by applying to each $C^{p}$-function $f_{n}$ the classical
Morse-Sard theorem \cite{Sard42}. $\hfill\Box$

\bigskip

\noindent As we see in the next section, several important classes
of lower semicontinuous functions satisfy the assumptions (thus
also the conclusions) of Proposition \ref{projec} and of Corollary
\ref{cor1}.

\section{Kurdyka-\L ojasiewicz inequalities for o-minimal functions}

\label{Kurdik}

Let us recall briefly a few definitions concerning o-minimal
structures (see for instance, Coste\,\cite{Coste99}, van\,der
Dries-Miller \cite{Dries-Miller96}, Ta\,L\^{e}\,Loi \cite{Taleloi98},
and references therein).

\begin{definition}
[o-minimal structure]\label{Definition_o-minimal}\rm{An o-minimal
structure on $(\mathbb{R},+,.)$ is a sequence of boolean algebras
$\mathcal{O}_{n}$ of ``definable'' subsets of $\mathbb{R}^{n}$,
such that for each $n\in\mathbb{N}$}
\smallskip

\rm{(i) if $A$ belongs to $\mathcal{O}_{n}$, then
$A\times\mathbb{R}$ and $\mathbb{R}\times A$ belong to
$\mathcal{O}_{n+1}$ ;}\smallskip

\rm{(ii) if $\Pi:\mathbb{R}^{n+1}\rightarrow\mathbb{R}^{n}$ is the
canonical projection onto $\mathbb{R}^{n}$ then for any $A$ in
$\mathcal{O} _{n+1}$, the set $\Pi(A)$ belongs to
$\mathcal{O}_{n}$ ;}\smallskip

\rm{(iii) $\mathcal{O}_{n}$ contains the family of algebraic
subsets of $\mathbb{R}^{n}$, that is, every set of the form
\[
\{x\in\mathbb{R}^{n}:p(x)=0\},
\]
where $p:\mathbb{R}^{n}\rightarrow\mathbb{R}$ is a polynomial
function ;}\smallskip

\rm{(iv) the elements of $\mathcal{O}_{1}$ are exactly the finite
unions of intervals and points. }
\end{definition}

\begin{definition}
[definable function]\label{Definition_def_fun}\rm{Given an
o-minimal structure $\mathcal{O}$ (over 
$(\mathbb{R},+,.)$), a function
$f:\mathbb{R}^{n}\rightarrow\mathbb{R}\cup\{+\infty\}$ is said to
be \emph{definable} in $\mathcal{O}$ if its graph belongs to
$\mathcal{O}_{n+1}$.}
\end{definition}

\begin{remark}
\rm{At a first sight, o-minimal structures might appear artificial
in optimization. The following properties (see
\cite{Dries-Miller96} for the details) might convince the reader
that this is not the case.} \smallskip

\rm{(i) The collection of \emph{semialgebraic sets} is an
o-minimal structure. Recall that semialgebraic sets are Boolean
 combinations of sets of the form
\[
\{x\in\mathbb{R}^{n}:p(x)=0,q_{1}(x)<0,\ldots,q_{m}(x)<0\},
\]
where $p$ and $q_{i}$'s are polynomial functions on
$\mathbb{R}^{n}$.}\smallskip

\rm{(ii) There exists an o-minimal structure that contains the
sets of the form}
\[
\{(x,t)\in\lbrack-1,1]^{n}\times\mathbb{R}:f(x)=t\}
\]
\rm{where $f:\mathbb{R}^{n}\rightarrow\mathbb{R}$ is real-analytic
around $[-1,1]^{n}$.} \smallskip

\rm{(iii) There exists an o-minimal structure that contains
simultaneously the graph of the exponential function
$\mathbb{R}\ni x\mapsto\exp x$ and all semialgebraic sets
(respectively, the structure defined in (ii)).} \smallskip

\rm{Let us finally recall the following important fact: the
composition of mappings that are definable in some o-minimal
structure remains in the same structure
\cite[Section\,2.1]{Dries-Miller96}. This is true for the sum, the
inf-convolution and several other classical operations of analysis
involving a finite number of definable objects. This remarkable
stability, combined with new techniques of finite-dimensional
optimization offers a large field of investigation. Several works
have already been developed in this spirit, see for instance
\cite{Drummond}, \cite{AMA2004}, \cite{BDLS2005}.}
\end{remark}

Given any o-minimal structure $\mathcal{O}$ and any lower
semicontinuous definable function
$f:\mathbb{R}^{n}\rightarrow\mathbb{R}\cup\{+\infty\}$ the
assumptions of Proposition \ref{projec} are satisfied. More
precisely, we have the following result.

\begin{lemma}
\label{Lemma_o-minimal_strat} Let
$\mathcal{B}:=\{B_{1},\ldots,B_{p}\}$ be a collection of definable
subsets of $\mathbb{R}^{n}$. Then there exists a definable
$C^{p}$-Whitney stratification $\{S_{1},\ldots,S_{\ell}\}$ of the
graph $\mbox{Graph\,}f$ of $f$ satisfying the transversality
condition $(\mathcal{H})$ and yielding $($by projecting each
stratum $S_{i}\subset \mathbb{R}^{n+1}$ onto $\mathbb{R}^{n})$ a
$C^{p}$-Whitney stratification $\{X_{1},\ldots,X_{\ell}\}$ of the
domain $\mathrm{dom\,}f\ $of $f$ compatible with $\mathcal{B}$.
\end{lemma}

\noindent\textbf{Proof} Let $\{\sigma_{1},\ldots,\sigma_{m}\}$ be
a $C^{1} $-stratification of the definable set $\mbox{Graph\,}f$
and set $\Sigma _{i}=\Pi(\sigma_{i})$ and $Y_{i,j}=\Sigma_{i}\cap
B_{j}$ for $1\leq i\leq m$ and $1\leq j\leq p$. Since the mapping
 $\Pi:\sigma_i\rightarrow\Sigma_i$ is an open continuous one to one mapping,
 its inverse $\Sigma_i\ni x\rightarrow (x,f(x))\in \sigma_i$ is continuous
 and thus so is $f|_{\Sigma_{i}}$. Since the
restriction $f|_{Y_{i,j}}$ of $f$ to each definable set $Y_{i,j}$
is continuous,
there exists a $C^{1}$-stratification $\{Z_{i,j,k}\}_{1\leq k\leq
q_{i.j}}$ of each definable set $Y_{i,j}$ such that
$f|_{Z_{i,j,k}}$ is $C^{1}$ (see \cite[Theorem 6.7]{Coste99},
\cite{Kaloshin2005} or \cite{Taleloi98} for example). Set
\[
A_{i,j,k}=\{(x,f(x)):x\in Z_{i,j,k}\}
\]
and consider a $C^{p}$-Whitney stratification $\mathcal{S}=\{S_{1}
,\ldots,S_{\ell}\}$ of the definable set $\mbox{Graph\,}f$
compatible with the definable sets $A_{i,j,k}$ (see \cite[Theorem
1.3]{Taleloi98}, for example). It is easily seen that
$\mathcal{S}$ satisfies the transversality condition
($\mathcal{H}$). Indeed, since
$\mathrm{dom\,}f=\bigcup\nolimits_{i,j,k} Z_{i,j,k}$, the relation
$e_{n+1}\in T_{(x,f(x))}S_{i}$ for some $x\in\mathrm{dom\,}f$
belonging say to the stratum $Z_{i,j,k}\subset S_{i}$ is
contradicting the differentiability of $f|_{Z_{i,j,k}}$ at $x$.
Setting $X_{i}=\Pi(S_{i})$, it is easily seen that the obtained
$C^{p}$-Whitney stratification
$\mathcal{X}=\{X_{1},\ldots,X_{\ell}\}$ of $\mathrm{dom\,} f$ is
compatible with the collection $\mathcal{B}$.$\hfill\Box$

\bigskip

\begin{remark}
\label{Remark_strat-fun} \rm{The aforementioned result can also be
obtained by evoking more delicate results on stratification of
functions (\cite{Taleloi97}, for example). We  give an
elementary proof for the reader's convenience.}
\end{remark}

\begin{corollary}
\label{Cstrat} Let
$f:\mathbb{R}^{n}\rightarrow\mathbb{R}\cup\{+\infty\}$ be a lower
semicontinuous definable function. There exists a finite definable
Whitney stratification $\mathcal{X}=(X_{i})_{i\in I}$ of
$\mbox{dom\,}f$ such that for all $x\in\mbox{dom\,}f$
\begin{equation}
\mathrm{Proj\,}_{T_{x}X_{x}}\mathrm{\,}\partial^{\circ}f(x)\mathrm{\,}
\subset\mathrm{\,}\{\mathrm{\,}\nabla_{\emph{R}}f(x)\mathrm{\,}
\}.\label{projfor}
\end{equation}
As a consequence\smallskip

(i) For all $x\in\mbox{dom\,}\partial^{\circ}f$ and $x^{\ast}\in
\partial^{\circ}f(x)$, we have $||\nabla_{\emph{R}}f(x)||\,\leq\,||x^{\ast
}||\,;$\smallskip

(ii) The set of Clarke critical values of $f$ is finite $;$\smallskip

(iii) The set of asymptotic Clarke critical values of $f$ is finite.
\end{corollary}

\noindent\textbf{Proof} Assertion (i) is a direct consequence of
(\ref{projfor}). This projection formula follows directly by
combining Lemma \ref{Lemma_o-minimal_strat} with Proposition
\ref{projec}. To prove (iii), let $f_{i}$ be the restriction of
$f$ to the stratum $X_{i}$. Then assertion (i), together with the
fact that the number of strata is finite, implies that the set of
the asymptotic Clarke critical values of $f$ is the union (over
the finite set $I$) of the asympotic critical values of each
(definable $C^{1}$) function $f_{i}$. Thus the result follows from
\cite[Remarque 3.1.5]{Dacunto}. Assertion (ii) follows directly
from (iii) (cf. Remark \ref{Remark1} (iii)).$\hfill\Box$

\bigskip

\begin{remark}
\rm{The fact that the set of the asymptotic critical values of a
definable differentiable function $f$ is finite has been
established in \cite[Th\'{e}or\`{e}me 3.1.4]{Dacunto} (see also
\cite[Theorem\,3.1]{KOS2000} for the case that the domain of $f$
is bounded). In \cite[Proposition 2]{Kurdyka98} a more general
result (concerning functions taking values in $\mathbb{R}^{k}$)
has been established in the semialgebraic case.}
\end{remark}

Before we proceed, let us recall from Kurdyka \cite[Theorem 1]{Kurdyka98} the
following result:

\begin{theorem}
[Kurdyka-\L ojasiewicz inequality]\label{Thm_KLI}Let $f:U\rightarrow
\mathbb{R}_{+}$ be a definable differentiable function, where $U$ is an open
and bounded subset of $\mathbb{R}^{n}$. Then there exist $\rho,c>0$ and a
strictly increasing definable function $\psi:(0,\rho)\rightarrow(0,+\infty)$
of class $C^{1}$ such that
\begin{equation}
||\nabla(\psi\circ f)(x)||\geq c,\text{\qquad for each }x\in U\cap
f^{-1}(0,\rho)\text{.}\label{Kurd}
\end{equation}
\end{theorem}

\begin{remark}
\label{rem-kurdyka} \rm{Let us observe that in the conclusion of
the above result, there is no loss of generality to assume $c=1$
and $\psi$ being defined and continuous on $[0,\rho)$ with
$\psi(0)=0.$ Moreover, a careful examination of the proof of
\cite[Theorem 1]{Kurdyka98} shows that the result of Theorem
\ref{Thm_KLI} remains valid if $U$ is any nonempty bounded
definable submanifold of $\mathbb{R}^{n}$.}
\end{remark}

We shall use Corollary \ref{Cstrat} to extend Theorem \ref{Thm_KLI} to a
nonsmooth setting.

\begin{theorem}
[Nonsmooth Kurdyka-\L ojasiewicz
inequality]\label{Theorem_KLineq}Let
$f:\mathbb{R}^{n}\rightarrow\mathbb{R\cup\{+\infty\}}$ be a lower
semicontinuous definable function and $U$ be a bounded definable
subset of $\mathbb{R}^{n}$. There exist $\rho>0$ and a strictly
increasing continuous definable function
$\psi:[0,\rho)\rightarrow(0,+\infty)$ which is $C^{1}$ on
$(0,\rho)$ with $\psi(0)=0$ and such that for all
$x\in\,U\,\cap\,|\,f\,|^{-1}(0,\rho)$ and all
$x^{\ast}\in\partial^{\circ}f(x)$
\begin{equation}
||x^{\ast}||\,\geq\,\frac{1}{\psi^{\prime}(|f(x)|)}.\label{KL}
\end{equation}
\end{theorem}

\noindent\textbf{Proof} Set
$U_{1}=\{x\in\,U\,\cap\,\mathrm{dom\,}f:f(x)>0\}$ and
$U_{2}=\{x\in\,U\,\cap\,\mathrm{dom\,}f:f(x)<0\}$ and let
$X_{1},\dots,X_{l}$ be a finite definable stratification of
$\mathrm{dom\,}f$ compatible with the bounded (definable) sets
$U_{1}$ and $U_{2}$ such that the definable sets
$S_{i}=\{(x,f(x):x\in X_{i}\}$ are the strata of a definable
$C^{p}$-Whitney stratification of $\mbox{Graph\,}f$ satisfying
($\mathcal{H}$) (cf. Lemma \ref{Lemma_o-minimal_strat}). For each
$i\in\{1,\ldots,l\}$ such that $X_{i}\subset U_{1}$ we consider
the positive $C^{1}$ function $f_{i} :=f|_{X_{i}}$ on the
definable manifold $X_{i}$ (thus for $x\in X_{i}$ we have $\nabla
f_{i}(x)=\nabla_{\emph{R}}f(x)$ and $f_{i}(x)=f(x)$) and we apply
Theorem \ref{Thm_KLI} (and Remark\,\ref{rem-kurdyka}) to obtain
$\rho_{i}>0$ and a strictly increasing definable $C^{1}$-function
$\psi_{i}:(0,\rho _{i})\rightarrow(0,+\infty)$ such that for all
$x\in f_{i}{}^{-1}(0,\rho_{i})$ we have
$||\nabla_{\emph{R}}f(x)||\geq\lbrack\psi_{i}^{\prime}(f(x))]^{-1}$.
Similarly, for each $j\in\{1,\ldots,l\}$ such that $X_{j}\subset
U_{2}$ we consider the positive $C^{1}$ function
$f_{j}:=-f|_{X_{i}}$ (note that for $x\in X_{j}$ we have $\nabla
f_{j}(x)=-\nabla_{\emph{R}}f(x)$ and $f_{j}(x)=-f(x)$) to obtain
as before $\rho_{j}>0$ and a strictly increasing definable
$C^{1}$-function $\psi_{j}:(0,\rho_{j})\rightarrow(0,+\infty)$
such that for all $x\in f_{j}{}^{-1}(0,\rho_{i})$ we have
$||\nabla_{\emph{R}
}f(x)||\geq\lbrack\psi_{j}^{\prime}(-f(x))]^{-1}$. Thus for all
$i\in \{1,\ldots,l\}$ there exist $\rho_{i}>0$ and a strictly
increasing definable $C^{1}$-function
$\psi_{i}:(0,\rho_{i})\rightarrow\mathbb{R}$ such that
\[
||\nabla_{\emph{R}}f(x)||\geq\frac{1}{\psi_{i}^{\prime}(|f(x)|)}
,\qquad\text{for all }x\in\,U\,\cap\,|f|{}^{-1}(0,\rho_{i}).
\]
Set $\rho=\min\rho_{i}$ and let $i_{1},i_{2}\in\{1,\ldots,l\}$. By the
monotonicity theorem for definable functions of one variable (see \cite[Lemma
2]{Kurdyka98}, for example), the definable function
\[
(0,\rho)\ni r\rightarrow1/\psi_{i_{1}}^{\prime}(r)-1/\psi_{i_{2}}^{\prime}(r)
\]
has a constant sign in a neighborhood of $0$. Repeating the argument for all
couples $i_{1},i_{2}$ and shrinking $\rho$ if necessary, we obtain the
existence of a strictly increasing, positive, definable function $\psi
=\psi_{i_{0}}$ on $(0,\rho)$ of class $C^{1}$ that satisfies $1/\psi^{\prime
}\leq1/\psi_{i}^{\prime}$ on $(0,\rho)$ for all $i\in\{1,\ldots,l\}$. Evoking
Corollary \ref{Cstrat}\thinspace(i), we obtain for all $x\in\,U\,\cap
\,|f|^{-1}(0,\rho)$ and all $x^{\ast}\in\partial^{\circ}f(x)$
\[
||x^{\ast}||\geq||\nabla_{\emph{R}}f(x)||\geq\frac{1}{\psi^{\prime}(|f(x))|}.
\]
Since $\psi$ is definable and bounded from below, it can be extended
continuously to $[0,\rho).$ By adding eventually a constant, we can also
assume $\psi(0)=0$.$\hfill\Box$

\bigskip

The assumption that the function $f$ is definable is important for the
validity of (\ref{KL}). It implies in particular that the connected components
of the set of the Clarke critical points of $f$ lie in the same level set of
$f$ (cf. Corollary \ref{Cstrat}\,(ii)). Let us present some examples of
$C^{1}$-functions for which (\ref{KL}) is not true.

\begin{example}
\label{Example_contrexemple} \rm{(i) Consider the function
$f:\mathbb{R} \rightarrow\mathbb{R}$ with}
\[
f(x)=\left\{
\begin{array}
[c]{cc}
x^{2}\sin\frac{1}{x}, & \text{ if }x\neq0\\
0, & \text{ if }x=0
\end{array}
\right.
\]
\rm{Then the set $S=\{x\in\mathbb{R}:f^{\prime}(x)=0\}$ meets
infinite many level sets. Consequently, (\ref{KL}) is not
fulfilled.}

\smallskip

\noindent\rm{(ii) A nontrivial example is proposed in \cite[page
14]{PalDem82}, where a $C^{\infty}$ ``Mexican-hat'' function has
been defined. An example of a similar nature has been given in
\cite{AMA2004}, and will be described below: Let $f$ be defined in
polar coordinate on $\mathbb{R}^{2}$ by}
\[
f(r,\theta)=\left\{
\begin{array}
[c]{cc}
\exp(-\frac{1}{1-r^2}\,)\,[1-\frac{4r^{4}}{4r^{4}+(1-r^{2})^{4}}\sin
(\theta-\frac{1}{1-r^{2}})], & \quad\text{if }r\leq1\\
& \\
0, & \quad\text{if }r>1.
\end{array}
\right.
\]
\rm{The function $f$ does not satisfy the Kurdyka-\L ojasiewicz
inequality for the critical value $0$, \textit{i.e.} one can not
find a strictly increasing $C^{1}$ function
$\psi:(0,\rho)\rightarrow(0,+\infty)$, with $\rho>0$, such that}
\[
||\nabla(\psi\circ f)(x)||\geq1
\]
\rm{for small positive values of $f(x)$. To see this, let us
notice that the proof of (\cite[Theorem 2]{Kurdyka98}) shows that
for \emph{any } C$^{1}$ function $f$ (not necessarily definable)
that satisfies the Kurkyka-\L ojasiewicz inequality, the bounded
trajectories of the gradient system}
\[
\dot{x}(t)+\nabla f(x(t))=0
\]
\rm{have a bounded length. However, in the present example, taking
as initial condition $r_{0}\in(0,1)$ and $\theta_{0}$ such that
$\theta_{0}(1-r_{0} )^{2}=1$, the gradient trajectory
$\dot{x}(t)=-\nabla f(x(t))$ must comply with}
\[
\theta(t)=\frac{1}{1-r(t)^{2}},
\]
\rm{where $r(t)\nearrow1^{-}$ as $t\rightarrow+\infty$ (see
\cite{AMA2004} for details). The total length of the above curve
is obviously infinite, which shows that the Kurdyka-\L ojasiewicz
inequality (for the critical value $0$) does not hold.}
\end{example}

\bigskip

Let us finally give an easy consequence of Theorem \ref{Theorem_KLineq} for
the case of subanalytic functions.

\begin{corollary}
[Subgradient inequality]\label{Corollary_Loja} Assume that
$f:\mathbb{R} ^{n}\rightarrow\mathbb{R\cup\{+\infty\}}$ is a lower
semicontinuous globally subanalytic function and $f(x_{0})=0$. Then there
exist $\delta>0$ and $\theta\in\lbrack0,1)$ such that for all
$x\in\,|\,f\,|^{-1}(0,\delta)$ we have
\[
|f(x)|^{\theta}\,\leq\,\rho\,||x^{\ast}||,\qquad\text{for all\
}x^{\ast} \in\partial^{\circ}f(x).
\]
\end{corollary}

\noindent\textbf{Proof} In case that $f$ is globally subanalytic, one can apply
\cite[Theorem\,(LI)]{Kurdyka98} to deduce that the continuous function $\psi$
of Theorem \ref{Theorem_KLineq} can be taken of the form $\psi(s)=s^{1-\theta
}$ with $\theta\in(0,1)$. $\hfill\Box$

\bigskip

\begin{remark}
\label{Remark_convex-stable}\rm{Corollary \ref{Cstrat}(ii) (and a
fortiori Corollary \ref{Corollary_Loja}) extends \cite[Theorem
7]{BDLS2005} to the lower semicontinuous case. We also remark that
the conclusions of Theorem \ref{Theorem_KLineq} and of Corollary
\ref{Corollary_Loja} remain valid for any notion of
subdifferential that is included in the Clarke subdifferential,
thus, in particular, in view of (\ref{inclusions}), for the
Fr\'{e}chet and the limiting subdifferential. However, let us point
out that this is not the case for broader notions of
subdifferentials, as for example the \emph{convex-stable}
subdifferential introduced and studied in \cite{BLO2002}. It is
known that the convex-stable subdifferential coincides with the
Clarke subdifferential whenever the function $f$ is locally
Lipschitz continuous, but it is strictly larger in general,
creating more critical points. In particular,
\cite[Section\,4]{BDLS2005} constructs an example of a subanalytic
continuous function on $\mathbb{R}^{3}$ that is strictly
increasing in a segment lying in the set of its broadly critical
points (that is, critical in the sense of the convex-stable
subdifferential). Consequently, Theorem \ref{Theorem_KLineq} and
Corollary \ref{Corollary_Loja} do not hold for this
subdifferential.}
\end{remark}

\bigskip

\noindent\textbf{Acknowledgment} The second author wishes to thank
K.\,Kurdyka and S.\,Simon for useful discussions. A part of this
work has been done during a visit of the first author at the
C.R.M. (Universitat Aut\`{o}noma de Barcelona). The first author
wishes to thank the C.R.M. for the financial support.

\begin{center}
\bigskip

----------------------------------------------------
\end{center}

\noindent J\'{e}r\^{o}me BOLTE \quad(\texttt{bolte@math.jussieu.fr}
;\texttt{\ http://www.ecp6.jussieu.fr/pageperso/bolte/})\smallskip

\noindent Equipe Combinatoire et Optimisation (UMR 7090), Case 189,
Universit\'{e} Pierre et Marie Curie\newline 4 Place Jussieu, 75252 Paris
Cedex 05.

\bigskip

\noindent Aris DANIILIDIS\quad(\texttt{arisd@mat.uab.es}\thinspace
;\quad\texttt{http://mat.uab.es/\symbol{126}arisd})\smallskip

\noindent Departament de Matem\`{a}tiques, C1/320\newline Universitat
Aut\`{o}noma de Barcelona\newline E-08193 Bellaterra (Cerdanyola del
Vall\`{e}s), Spain.

\bigskip

\noindent Adrian LEWIS\quad(\texttt{aslewis@orie.cornell.edu}\thinspace
;\quad\texttt{http://www.orie.cornell.edu/\symbol{126}aslewis})\smallskip

\noindent School of Operations Research and Industrial Engineering\newline
Cornell University\newline 234 Rhodes Hall, Ithaca, NY 14853, United States.

\bigskip

\noindent Masahiro SHIOTA\quad(\texttt{shiota@math.nagoya-u.ac.jp}) \smallskip

\noindent Department of Mathematics\newline Nagoya University\, (Furocho,
Chikusa) \newline Nagoya 464-8602, Japan.
\end{document}